\documentclass[reqno,draft]{amsart}
\usepackage{multicol, color}

\newcommand{\rem}[1]{}

\theoremstyle{plain}
\newtheorem{lemma}{Lemma}
\newtheorem{theorem}[lemma]{Theorem}
\newtheorem{corollary}[lemma]{Corollary}

\newtheorem{definition}[lemma]{Definition}
\theoremstyle{remark}
\newtheorem{remark}{Remark}

\newcommand*  {\N} {{\mathbb N}}

\newcommand*  {\R} {{\mathbb R}}

\def\Om{\Omega}

\def\pp{\partial}
\def\RN{\R^N}
\def\vpr{\varphi_R}
\def\vp{\varphi}
\def\gd{\nabla}
\def\lp{\triangle}
\def\bu{\bold{u}}

\begin{document}
\title[steady state of Kuramoto-Sivashinsky equations]
{trivial stationary solutions to the Kuramoto-Sivashinsky and certain nonlinear
  elliptic equations}
\date{March 30, 2006}

\author[Y.Cao]{Yanping Cao}
\address[Y.Cao]
{Department of Mathematics\\
University of California\\
Irvine, CA 92697-3875,USA}
\email{ycao@math.uci.edu}

\author[E.S. Titi]{Edriss S. Titi}
\address[E.S. Titi]
{Department of Mathematics \\
and  Department of Mechanical and  Aerospace Engineering \\
University of California \\
Irvine, CA  92697-3875, USA \\
{\bf ALSO}  \\
Department of Computer Science and Applied Mathematics \\
Weizmann Institute of Science  \\
Rehovot 76100, Israel}
\email{etiti@math.uci.edu and edriss.titi@weizmann.ac.il}

\begin{abstract}
We show that the only locally integrable stationary solutions to the integrated Kuramoto-Sivashinsky
equation in $ \R$  and $ \R^2$  are the trivial constant solutions.  We extend our technique and prove
similar results to other
nonlinear elliptic problems in $\R^N$.
\end{abstract}

\maketitle



\section{Introduction}   \label{S-1}
The integrated Kuramoto-Sivashinsky equation (abbreviated hereafter as the KSE)
\begin{equation}
\hskip-.8in
\phi_t+\lp^2\phi+\lp\phi+\frac{1}{2}|\gd\phi |^2=0   \label{KSE-1}
\end{equation}
subject to appropriate initial and boundary conditions has been
introduced in \cite{Kuramoto},\cite{Kuramoto-Tsuzuki} and in
\cite{Sivashinskyintro1},\cite{Sivashinskyintro2} in studying phase
turbulence and the flame front propagation in combustion theory. In
the absence of any a priori estimates for the solutions of the
scalar equation (\ref{KSE-1}), most authors find it more convenient,
for the mathematical study, to consider the differential form of the
equation for $\bu=\gd \phi$
\begin{equation}
\hskip-.8in
\bu_t+\lp^2\bu+\lp\bu+(\bu\cdot\gd)\bu =0.      \label{KSE-2}
\end{equation}

The one-dimensional case has been studied by many authors (see e.g.,
\cite{Bronski}, \cite{CheskidovFoias}, \cite{Collet}, \cite{CFNT},
\cite{FK}, \cite{FNST}, \cite{FST}, \cite{Goodman},
\cite{Ilyashenko}, \cite{JKT}, \cite{KNS}, \cite{DMichelson},
\cite{NST}, \cite{Tadmor}, \cite{Temam} and references therein). It
is known that the long-term dynamics of the 1-D equation
(\ref{KSE-2}) with periodic boundary condition, of period $L$,
possesses a  complicated global attractor $\mathcal{A}_L$ with
finite dimension (see e.g., \cite{CFNT}, \cite{Goodman},
\cite{Ilyashenko}, \cite{DMichelson}, \cite{NST}, \cite{Temam} and
references therein). The best upper bound for the dimension of the
global attractor is of the order  $0(L^\frac{45}{40})$ is obtained
based on the best available upper bound for the size of the
absorbing ball \cite{FOtto}. Namely, the current estimates for the
upper bound for the dimension of the global attractor depend
explicitly on $R=\limsup_{t \rightarrow \infty}|u(\cdot,t)|_{L^2}$:
if $R$ is of the order $0(L^\beta)$, then the upper bound for the
Hausdorff and fractal dimensions
 of the global attractor satisfies $d_H(\mathcal{A}_L) \le d_f(\mathcal{A}_L) \le
0(L^\frac{30+10\beta}{40})$. As mentioned above, the best estimate
for $R$ is given in \cite{FOtto}: $R \sim o(L^\frac{3}{2})$ (see
also \cite{Bronski}). On the other hand, based on numerical
simulations and physical arguments it is conjectured \cite{Pomeau}
that the upper bound for the dimension of the global attractor
should behave like $L$. This conjectured estimate also matches the
readily available lower bound for the dimension of the global
attractor which is obtained by linearizing about the stationary
solution $u\equiv 0$. To achieve this conjectured bound for the
dimension of the global attracto it requires, for instance, to
establish a uniform bound for the $L^\infty$-norm of the solutions
on the attractor, which is independent of the period $L$. Therefore,
the question is: whether the $L^\infty$-norms of any solution
$u(x,t)$ on the attractor are uniformly bounded, independent of $L$.
The remarkable paper of Michelson \cite{DMichelson} shows that all
the bounded stationary solutions to equation (\ref{KSE-2}) on the
whole line are uniformly bounded by a constant $K_M$. In a paper by
Cheskidov and Foias \cite{CheskidovFoias}, they consider the
nonhomogeneous one-dimensional stationary KSE (\ref{KSE-2}) subject
to periodic boundary condition with zero spatial average. Namely,
the problem
\begin{eqnarray*}
&&\hskip-.8in
u_{xxxx}+u_{xx}+uu_x=f(x)                   \\
&&\hskip-.8in
u(x)=u(x+L), \quad   \int^L_0 u(x)dx=0.
\end{eqnarray*}
They find an explicit estimate for the Michelson constant $K_M$,
namely, $K_M\le92.2$ . Furthermore, they study the set of averages
of the solutions with respect to invariant probability time average
measures:
\begin{equation*}
\hskip-.8in E=\left\{{\bar{u}}=\int_H u\mu (du)\| \mu \mbox{ is
invariant probability time average measure on }\mathcal{A}_L,L>0
\right\},
\end{equation*}
where $H=\{u: u \in L^2,u(x)=u(x+L),\int_0^Lu(x)dx=0\}$ the phase
space. For details about stationary statistical solutions and
invariant probability time average measures see, e.g., \cite{CFMRT}.
In particular, it is  shown in \cite{CheskidovFoias} that the
$L^\infty$-norm of the set E, defined above, is uniformly bounded,
independent of L. Indeed, Michelson's result is a particular case of
the above set. This is because for $\mu=\delta_{u_0}$, when $u_0$ is
a stationary solution to the KSE, gives an invariant probability
time average measure. In this case,
\begin{equation*}
\hskip-.8in
\bar{u}=\int_H u\:\mu (du)=\int_H u\: \delta_{u_0} (du)=u_0.
\end{equation*}

The question of global regularity of the Cauchy problem
\begin{eqnarray}
&&\hskip-.8in
\phi_t + \lp ^2 \phi + \gd \phi+ \frac{1}{2} |\gd \phi |^2 =0  \qquad \mbox {in}\quad \RN \nonumber\\
&&\hskip-.8in
\phi(x,0)=\phi_0(x)                                    \label{Cauchy}
\end{eqnarray}
or the periodic boundary condition case
\begin{eqnarray}
&&\hskip-.8in
\phi_t+ \lp ^2 \phi +\gd \phi + \frac{1}{2} |\gd \phi |^2=0    \nonumber \\
&&\hskip-.8in
\phi(x+Le_j,t)=\phi(x,t)   \qquad \qquad \mbox{for}\quad j=1,2,...,N                          \label{periodic}\\
&&\hskip-.8in
\phi(x,0)=\phi_0(x)                                                              \nonumber
\end{eqnarray}
is still an open question in dimensions two and higher cases (see,
however, \cite{SellTaboada} for the case  of thin two dimensional
domains for large, but restricted, initial data). Motivated by the
question of global regularity of (\ref{Cauchy}) or (\ref{periodic}),
the authors of \cite{Bellout} study what they call the hyper-viscous
Hamiltion-Jacobian-like equation for the scalar function u:
\begin{eqnarray}
&&\hskip-.8in
u_t+\lp ^2 u=|\gd u|^p                         \nonumber \\
&&\hskip-.8in
u|_{\pp\Om}=\gd u|_{\pp\Om}=0                  \label{HJ}\\
&&\hskip-.8in
u(x,0)=u_0(x),                                  \nonumber
\end{eqnarray}
where $\Om$ is a smooth bounded domain  in $\RN$.

In the case $p> 2$ , they show that finite time blow-up will occur
for special ``large" initial condition. It is remarked that the
blow-up occurs in $L^\infty$-norm, i.e., the derivative of the u
remains finite as long as the solution exists and has finite
$L^\infty$-norm. In particular, there is essential difference in the
structure of formation of singularities from that of generalized
viscous Hamiltion-Jacobi equations \cite{Souplet}:
\begin{eqnarray}
&&\hskip-.8in
u_t-\lp u=|\gd u|^p    \qquad  \mbox{ in } \quad  \Om \times (0,\infty ) \nonumber \\
&&\hskip-.8in
u|_{\partial \Om}=0                                 \label{GHJ}\\
&&\hskip-.8in
u(x,0)=u_0(x)         \quad\quad\quad \mbox{   in   } \quad\Om ,                     \nonumber
\end{eqnarray}
where $\Om$ is a smooth bounded domain in $\RN$. Regradless of the
value of $p$, $p\geq 0$,  problem (\ref{GHJ}) satisfies a maximum
principle, the $L^\infty$-norm of the solutions to problem
(\ref{GHJ}) remains bounded for as long as the solutions exist.
Thus, the solutions to (\ref{GHJ}) that become singular in finite
time must develop their singularities in one of their derivatives
(see \cite{Souplet}). However, for the critical case of $p=2$ in
problem (\ref{HJ}), it is still unknown whether there is global
regularity or there will be finite time blow-up for certain initial
data. It is worth mentioning, however, that for the case $p=2$ the
problem (\ref{HJ}) is the viscous Burgers equation which globally
well-posed for $N=1,2,3.$

Motivated by the above discussion, we study in this paper the steady
state problem (\ref{KSE-1}) in $\R^N$ for $N\ge 1$. In particular,
we show that in dimensions $N=1,2$ the only locally integrable
steady solutions of (\ref{Cauchy}) or (\ref{periodic}) are the
trivial solutions $\phi(x)=\mbox{constant}$. The techniques
developed and used here are inspired by the work of Mitidieri and
Pohozaev \cite{Pohozaev}. It is worth mentioning that for $N=3$,
Michelson \cite{Michelson} has established, using asymptotic
methods, the existence of a nontrivial radial steady state solution
of (\ref{Cauchy}). This is consistent with our results which are
restricted to dimensions $N=1,2$.

In section (\ref{SEC-2}), we will study the stationary solutions and present our main result.
In section (\ref{SEC-3}), we extend our tools to certain nonlinear elliptic systems in $\RN$
and show nonexistence of nontrivial solutions to those equations.

\section{steady state Kuramoto-Sivashinsky equation}   \label{SEC-2}
In this section we consider the integrated version of the homogeneous steady state KSE in $\RN$
\begin{equation}
\lp^2\phi+\lp\phi+\frac{1}{2}|\gd\phi|^2=0.  \label{SSKSE-1}
\end{equation}

We emphasize the fact that we do not require $\phi$ to satisfy any specific boundary condition as $|x|\rightarrow \infty$.
\begin{definition}
A function $\phi \in H^1_{loc}(\RN)$ is  called a locally integrable solution of (\ref{SSKSE-1}) in $\RN$ if $\phi$ satisfies
the equation (\ref {SSKSE-1}) in the distribution sense, i.e., in $\mathcal{D'}(\RN)$.
\end{definition}

\begin{theorem}
For $N=1,2$, the only locally integrable solutions of equation (\ref {SSKSE-1}) are the trivial solution, i.e., $\phi$=constant.
\end{theorem}

\begin{proof}
Consider the smooth radial cut-off function $\varphi_0(x)\in C^\infty_0 (\RN), 0\le \varphi_0(x)\le 1$ such that
\begin{equation*}
\vp_0(x)=\left \{ \begin{array}{lll}
                 1 & |x|\le 1 \\
                 0 & |x|\ge 2 \\
                 \mbox{smooth} & 1<|x|<2.
                 \end{array}
         \right.
\end{equation*}
Let
\begin{equation}
\hskip-.8in
\vpr(x)=\varphi_0(\frac{x}{R})  \label {test-fct}
\end{equation}
Suppose $\phi$ is a locally integrable solution of (\ref{SSKSE-1}),
taking action of (\ref{SSKSE-1}) on the test function $\vpr$ , we have

\begin{equation}
\hskip-.8in
\int_{\RN} |\gd \phi(x)|^2 \vpr(x) dx=-2\langle \lp^2\phi, \vpr \rangle -2 \langle \lp\phi, \vpr \rangle.     \label{proof-1}
\end{equation}
We estimate the right hand side of the above equality,
\begin{eqnarray*}
2\left | \langle \lp^2\phi, \vpr \rangle \right |
&=& 2 \left | \int_{\RN} -\gd \phi(x) \cdot \gd (\lp \vpr(x)) \: dx \right |  \\
&\le& 2 \int_{\RN} |\gd \phi(x)| |D^3\vpr(x)| \: dx \\
&\le& 2 \left( \int_{\RN} |\gd\phi(x)|^2 \vpr(x) \: dx \right)^\frac{1}{2} \left(\int_{\RN} \frac{|D^3 \vpr(x)|^2}{\vpr(x)} \:dx\right)^\frac{1}{2}
\end{eqnarray*}
where $D^k$ denotes a generic expression of the form
\begin{equation*}
\hskip-.8in
D^k u=\sum_ {|\alpha |=k}a_{\alpha} \frac{\partial^{|\alpha |} \:u}{\partial x_1^{\alpha_1} \: \partial x_2^{\alpha_2} \cdot \cdot \cdot \partial x_N^{\alpha_N}}
\end{equation*}
where $\alpha=\left(\alpha_1,\alpha_2,...,\alpha_N \right)$ is a multi-index and $a_\alpha$ are constants.

\noindent
Also we have
\begin{eqnarray*}
2\left | \langle \lp \phi, \vpr \rangle \right |
&=& 2\left | -\int_{\RN}\gd\phi(x)\cdot \gd \vpr(x) \: dx \right | \\
& \le& 2 \left( \int_{\RN} |\gd \phi(x)| | D \varphi_R (x)| \: dx \right) \\
&\leq& 2\left( \int_{\RN}|\gd \phi(x)|^2 \vpr(x) \: dx \right)^\frac{1}{2} \left( \int_{\RN}\frac{|\gd\vpr(x)|^2}{\vpr(x)}\: dx\right)^\frac{1}{2}
\end{eqnarray*}
The above estimates and (\ref{proof-1}) imply
\begin{eqnarray}
\int_{\RN}|\gd\phi(x)|^2 \vpr(x) \: dx
&\le& 2 \left( \int_{\RN} |\gd \phi(x)|^2 \vpr(x) \: dx \right)^\frac{1}{2} \left( \int_{\RN} \frac{ |D^3 \vpr(x)|^2}{\vpr(x)}\: dx \right)^\frac{1}{2}+\nonumber\\
&+& 2\left( \int_{\RN}|\gd \phi(x)|^2 \vpr(x) \: dx \right)^\frac{1}{2} \left( \int_{\RN}\frac{|\gd\vpr(x)|^2}{\vpr(x)}\: dx\right)^\frac{1}{2}   \label{proof-2}
\end{eqnarray}
By Young inequality, we reach
\begin{equation*}
\hskip-.8in
\int_{\RN}|\gd\phi(x)|^2\vpr(x) \: dx \leq 8 \left( \int_{\RN}\frac{|D^3\vpr(x) |^2}{\vpr(x)} \: dx +\int_{\RN}\frac{|D\vpr(x)|^2}{\vpr(x)} \: dx \right).
\end{equation*}
By our definition  $\varphi_R(x)=\varphi_0(\frac{x}{R})$. Let us change the variables  $x=R\xi$, then we obtain
\begin{equation}
\hskip-.8in
\int_{\RN}|\gd\phi(x)|^2\vpr(x) \: dx\leq 8 \left(R^{N-6}\int_{1<|\xi| < 2} \frac{|D^3\varphi_0(\xi)|^2}{\varphi_0(\xi)} \: d \xi
  +R^{N-2} \int_{1< |\xi| <2} \frac{|D\varphi_0(\xi)|^2}{\varphi_0(\xi)} \: d \xi\right)                  \label{proof-3}
\end{equation}
\vskip0.1in
Now, we further specialize in the choice of the test function  $\varphi_0$  such that the integrals on the right hand side of (\ref{proof-3}) are finite.
Then (\ref{proof-3}) implies
\begin{equation}
\hskip-.8in
\int_{\RN}|\gd\phi(x)|^2\vpr(x)\: dx\leq 8C_0 R^{N-6}+8C_1R^{N-2} \label{proof-4}
\end{equation}
where
\begin{equation}
\hskip-.8in
C_0=\int_{1<|\xi |<2} \frac{|D^3 \varphi_0(\xi)|}{\varphi_0(\xi)}\: d \xi \mbox{ ,   }
C_1=\int_{1<|\xi |<2} \frac{| D\varphi_0(\xi)|}{\varphi_0(\xi)} \: d\xi  \label{test-fct2}
\end{equation}
\vskip0.1in
\noindent
{\bf{The case $N=1$}}

Let us first consider $N=1$. Let  $r>0$  be fixed large enough, we
consider $R$  to be large enough such that  $R>4r$. From
(\ref{proof-4}) we conclude that
\begin{eqnarray}
\int_{|x|<r} \left|\frac{ d \phi(x)}{dx}\right|^2 \:dx
&\le& \int_{\R} \left|\frac{d \phi(x)}{dx}\right|^2 \varphi_R(x) \: dx   \nonumber  \\
&\le& 8C_0 R^{-5}+8C_1 R^{-1}.                                        \label{proof-ad}
\end{eqnarray}
Passing to the limit as $R\rightarrow\infty$ in (\ref{proof-ad}), we obtain that
\begin{equation*}
\hskip-.8in
\int_{|x|<r}|\frac{d \phi(x)}{dx}|^2\: dx=0
\end{equation*}
for every  $r>0$. Therefore $\frac{d\phi}{dx}=0$ and the assertion of the theorem is proved for $N=1$.


\vskip0.1in
{\bf{The case $N=2$}}

Now we consider the case  $N=2$.  In this case, the relation (\ref{proof-3}) implies that
\begin{equation*}
\hskip-.8in
\int_{\R^2}|\gd \phi(x)|^2 \vpr(x) \: dx \leq 8C_0R^{-4}+8C_1.
\end{equation*}
Choose as before $r>0$  fixed large enough, and let  $R>4r$. Frome the above we get
\begin{eqnarray*}
\hskip-.8in
\int_{|x|<r} |\gd \phi(x)|^2 \:dx
&\le& \int_{\R^2} |\gd \phi(x)|^2 \varphi_R(x) \:dx  \\
&\le& 8C_0R^{-4}+8C_1.
\end{eqnarray*}
Passing to the limit as  $R\rightarrow \infty$ ,  we obtain
\begin{equation}
\hskip-.8in
\int_{|x|<r} |\gd \phi(x) |^2 \: dx \leq 8C_1   \label{proof-5}
\end{equation}
for every  $r>0$.
By the Lebesgue monotone convergence theorem, we conclude that
\begin{equation}
\hskip-.8in
\gd \phi \in L^2\left( \R^2 \right)\label{constant1}
\end{equation}
and
\begin{equation}
\hskip-.8in
\int_{\R^2} | \gd \phi(x)|^2 \: dx \le 8C_1 \label {constant2}
\end{equation}
Now, let us return to inequality (\ref{proof-1}). Note that
\begin{equation}
\hskip-.8in
\mbox{supp}\{ D\vpr\}\subseteq \{ x\in \RN | R \le |x| \le 2R\}.  \label{proof-6}
\end{equation}
We estimate the right hand side of the relation (\ref{proof-1}),
\begin{eqnarray*}
2\left |\langle \lp ^2 \phi , \varphi_R \rangle \right|
&=& 2\left| \int_{\R^2}-\gd \phi(x) \cdot \gd(\lp \varphi_R(x)) \: dx\right| \\
&\le& 2 \int_{R<|x|<2R}  |\gd \phi(x)|  |D^3 \varphi_R(x) | \: dx \\
&\le& 2 \left( \int_{R<|x|<2R} |\gd \phi(x) |^2 \varphi_R(x) \:dx \right)^{\frac{1}{2}}
\left( \int_{\R^2} \frac{|D^3 \varphi_R(x) |^2}{\varphi_R(x)} \:dx \right)^{\frac{1}{2}}\\
&\le& 2R^{-2}\left( \int_{R<|x|<2R} |\gd \phi(x)|^2 \: dx \right) ^{\frac{1}{2}} \left(\int_{1<|\xi|<2} \frac{|D^3 \varphi_0(\xi)|^2}{\varphi_0(\xi)}\: d\xi \right)^{\frac{1}{2}}\\
&\le& 2C_0^{\frac{1}{2}}R^{-2} \left(\int_{R<|x|<2R} |\gd \phi(x) |^2\: dx \right)^{\frac{1}{2}},
\end{eqnarray*}
where $C_0$ is given in (\ref{test-fct2}) and in the above we changed the variable $x=R\xi$ and applied (\ref{constant1}).
Similary, for the other integral on the right hand side of (\ref{proof-1}),
\begin{eqnarray*}
2\left| \langle \lp \phi, \varphi_R \rangle \right|
&=&2\left| \int_{\R^2} -\gd \phi(x) \cdot \gd \varphi_R(x) \: dx  \right| \\
&\le& 2 \int_{R<|x|<2R} | \gd \phi(x) | |D \varphi_R(x) | \:dx  \\
&\le& 2 \left( \int_{R<|x|<2R}|\gd \phi(x) |^2 \varphi_R(x) \: dx \right)^{\frac{1}{2}} \left(\int_{\R^2} \frac{|D\varphi_R(x)|^2}{\varphi_R(x)}\:dx\right)^{\frac{1}{2}}\\
&\le& 2 \left(\int_{R<|x|<2R} |\gd \phi(x)|^2 \: dx \right)^{\frac{1}{2}} \left( \int_{1<|\xi|<2} \frac{|D\varphi_0(x)|^2}{\varphi_0(x)}\:d\xi \right)^{\frac{1}{2}}\\
&\le& 2C_1^{\frac{1}{2}} \left(\int_{R<|x|<2R} |\gd \phi(x)|^2 \: dx \right)^{\frac{1}{2}}
\end{eqnarray*}
where $C_1$ is given in (\ref{test-fct2}) and in the above we changed the variable  $x=R\xi$  and applied (\ref{constant2}).\\
\noindent
These estimates and (\ref{proof-1}) imply
\begin{eqnarray}
\hskip-.8in
\int_{\R^2}|\gd \phi(x)|^2 \varphi_R(x) \: dx
&\le & 2C_0^{\frac{1}{2}}R^{-4} \left(\int_{R<|x|<2R} |\gd \phi(x)|^2 \:dx \right)^{\frac{1}{2}} \nonumber\\
&+& 2C_1^{\frac{1}{2}} \left( \int_{R<|x|<2R} |\gd \phi(x)|^2 \:dx \right)^{\frac{1}{2}}    \label{proof-detail}
\end{eqnarray}
Passing to the limit as  $R\rightarrow \infty$  in (\ref{proof-detail}), by (\ref{proof-5}), (\ref{constant2}), and the Lebesgue Dominated
Convergence Theorem we obtain that
\begin{equation*}
\hskip-.8in
\int_{\R^2} |\gd\phi|^2 \: dx =0
\end{equation*}
Hence, $\phi(x)=$  constant.
\end{proof}

As a consequence of the above Theorem we have the following corollary:
\begin{corollary}
The only solutions to (\ref{periodic}), i.e., the only periodic solutions of equation (\ref{SSKSE-1}), are the constants.
\end{corollary}
Actually, one can prove this corollary in a direct trivial way and
for all $N$. In this case the set of test functions  $\mathcal{V}$
consists of all trigonometric polynomials. The function
$\varphi(x)=1 \in \mathcal{V}$ is a test function. Taking the action
of (\ref{SSKSE-1}) on $\varphi$ in the region $\Omega=[0,L]^N$ will
give us
\begin{equation*}
\hskip-.8in
\langle \lp^2 \phi, 1 \rangle + \langle \lp \phi, 1 \rangle +\frac{1}{2} \int_{\Omega} |\gd \phi(x) |^2 \: dx =0
\end{equation*}
By the peridocity of $\phi$ in $\Omega$, we obtain
\begin{equation*}
\hskip-.8in
\int_\Omega |\gd \phi|^2 \: dx = 0
\end{equation*}
It follows readily that $\gd \phi=0$, which implies that $\phi=\mbox{constant}$.\quad\quad\quad\quad\quad\quad\quad\quad\quad\quad\quad\quad\quad\quad\quad\quad\quad\qedsymbol

\vskip.1in
Next we consider the nonhomogeneous steady state:
\begin{equation}
\hskip-.8in
\lp^2 \phi + \lp \phi +\frac{1}{2} |\gd \phi |^2=f  \label{SSKSE-2}
\end{equation}
where $f(x)\in L^1_{loc}\left(\RN \right)$.
\begin{corollary}
Let $N\leq 2$, and $\liminf_{R\rightarrow \infty} \int_{\RN} f(x) \varphi_R(x) \: dx \le \mu_0<0$ for some
constant $ \mu_0$. Here $\varphi_R$ is specified in the manner of (\ref{test-fct}) and (\ref{test-fct2}). Then equation
(\ref{SSKSE-2}) has no locally integrable solutions, i.e., no solutions in $ H^1_{loc} (\RN)$.
\end{corollary}

\begin{proof}
Taking action of (\ref{SSKSE-2}) on the test function $\vpr$ defined in (\ref {test-fct}) and (\ref{test-fct2}), we get
\begin{equation*}
\hskip-.8in \int_{\RN} |\gd \phi(x) |^2 \vpr(x) \: dx=\int_{\RN}
f(x) \vpr(x) \: dx -2 \langle \lp ^2 \phi, \vpr \rangle -2\langle
\lp \phi, \vpr \rangle .
\end{equation*}
By (\ref{proof-2}), we reach
\begin{equation*}
\hskip-.8in
\int_{\RN} |\gd \phi(x) |^2 \vpr(x) \: dx \le 2 \int_{\RN} f(x) \vpr(x) \: dx +8 \left( \int_{\RN} \frac{| D^3 \vpr(x) |^2}{ \vpr(x) } \: dx + \int_{\RN} \frac{ |D \vpr(x)|} {\vpr(x)} \: dx \right )
\end{equation*}
Since $ \liminf_{R\rightarrow \infty} \int_{\RN} f(x) \vpr(x) \: dx \le \mu_0$, using the same argument
as in the proof of Theorem 2, we obtain
\begin{equation*}
\hskip-.8in \int_{\RN} | \gd \phi(x) |^2 \: dx  \le  2 \mu_0 < 0 ,
\end{equation*}
for $N=1,2$, which implies that we do not have any locally
integrable solutions for equation (\ref{SSKSE-2}).
\end{proof}

\section{Generalization to other Nonlinear elliptic problem} \label{SEC-3}
In this section we generalize the tools developed in the previous section and apply them to
certain class of nonlinear elliptic problems. Consider the nonlinear elliptic equation
\begin{equation}
\hskip-.8in \left(-\lp \right)^m u \pm |\gd^l \lp^n u|^p=0
\label{NE-1}
\end{equation}
defined in the whole space $ \RN$,  where $l=0$ or $l=1$, $2n+l \ge 0$, $2m>2n+l$ and $p>1$.
\begin{definition}
A function $u \in W^{2n+l, p}_{loc}(\RN)$ is called locally integrable solution of equation (\ref{NE-1}) if u satisfies equation
(\ref{NE-1}) in the distribution sense.
\end{definition}
We emphasize again that we do not require the solution $u$ to satisfy any specific boundary condition as $|x| \rightarrow \infty$.

\begin{remark}
A solution u of equation (\ref{NE-1}) is said to be trivial if $\gd^l\lp^n u=0$, and it is called nontrivial otherwise.
\end{remark}

\begin{theorem}
Let $N\leq \frac{(2m-(2n+l))p}{p-1}$, then the only locally integrable solutions of the equation (\ref{NE-1}) are the trivial
solutions.
\end{theorem}

\begin{proof}
Taking the action of  (\ref{NE-1}) on the test function $\vpr$ definde in (\ref {test-fct}), we have
\begin{equation*}
\hskip-.8in \langle (-\lp )^m u, \vpr \rangle \pm
\int_{\RN}|\gd^l\lp^n u(x)|^p \vpr(x) \: dx =0
\end{equation*}
Or
\begin{equation*}
\hskip-.8in \int_{\RN} |\gd^l\lp^n u(x)|^p \vpr (x) \: dx=\mp
\langle (- \lp )^m u, \vpr \rangle
\end{equation*}
By definition of distribution, we have
\begin{eqnarray}
\left | \int_{\RN}|\gd^l\lp^n u(x) |^p \vpr(x)  \: dx \right |
& =& \left |\mp\langle (-\lp )^m u, \vpr \rangle \right |        \nonumber \\
&=&  \left | \int_{\RN} \gd^l\lp^n u(x)\cdot D^{2m-(2n+l)}\vpr(x) \: dx \right | \nonumber\\
&\le & \left( \int_{\RN} |\gd^l\lp^n u(x)|^p \vpr(x) \: dx \right) ^\frac{1}{p} \left( \int_{\RN} \frac{|D^{2m-(2n+l)} \vpr(x) |^{p'}}{\vpr(x)^{p'-1}} \right) ^{\frac{1}{p'}} \label{NE-2}
\end{eqnarray}
where $p'$ is  the conjugate of $p$: $\frac{1}{p}+\frac{1}{p'}=1$.
So, we have
\begin{equation*}
\hskip-.8in
\left( \int_{\RN}|\gd^l\lp^nu(x)|^p \vpr(x) \: dx \right) \leq \left( \int_{\RN} \frac{|D^{2m-(2n+l)}\vpr(x) |^{p'}}{\left(\vpr(x)\right)^{p'-1}} \: dx \right)
\end{equation*}
Again, we change the variables $x=R\xi$. If we further specify $\varphi_0$ such that
\begin{equation}
\tilde{C}_0=\int_{1<|\xi|<2} \frac{| D^{2m-(2n+l)} \varphi_0(\xi)|^{p'}}{(\varphi_0(\xi))^{p'-1}} \:d\xi < \infty  \label{ctilde}
\end{equation}
We will have
\begin{equation}
\hskip-.8in
\left( \int_{\RN}|\gd^l\lp^nu(x)|^p \vpr \: dx \right) \leq \tilde{C}_0 \: R^\theta   \label{NE-3}
\end{equation}
where $\theta=N-(2m-(2n+l))p'$.
\vskip.1in
\noindent
{\bf{The case $N< \frac{(2m-(2n+l))p}{p-1}$}}

Let us first consider the case $N<\frac{(2m-(2n+l))p}{p-1}$, i.e., $\theta<0$. Let $r>0$ be fixed large
enough. We consider $R$ be large enough such that $R>4r$. From (\ref{NE-3}) we conclude that
\begin{eqnarray}
\hskip-.8in
\int_{|x|<r}| \gd^l\lp^n u(x) |^p \: dx
&\le & \int_{\RN} |\gd^l\lp^n u(x) |^p \vpr(x) \: dx    \nonumber \\
&\le & \tilde{C}_0 R^\theta                      \label {NE-4}
\end{eqnarray}
Passing to the limit as $R \rightarrow \infty$ in (\ref{NE-4}), we obtain that
\begin{equation*}
\hskip-.8in
\int_{|x|<r}| \gd^l\lp^n u(x) |^p \: dx \le 0
\end{equation*}
for every $r>0$.
By Lebesgue Momotone Convergence Theorem, we conclude that
\begin{equation*}
\hskip-.8in
\int_{\RN} |\gd^l\lp^n u(x) |^p \: dx =0
\end{equation*}
Then the assertion of the theorem is proved for $N< \frac{(2m-(2n+l))p}{p-1}$.
\vskip.1in
{\bf{The case $N=\frac{(2m-(2n+l))p}{p-1}$}}

Next, we consider $\N=\frac{(2m-(2n+l))p}{p-1}$, i.e., $\theta=0$. Then the relation (\ref{NE-3}) implies
that
\begin{equation*}
\hskip-.8in
\int_{\RN}|\gd^l\lp^n u(x)|^p \vpr(x) \: dx \le \tilde{C}_0
\end{equation*}
Choose as before $r>0$ fixed large enough and let $R>4r$. From the above relation we obtain
\begin{eqnarray*}
\hskip-.8in
\int_{|x|<r} |\gd^l\lp^n u(x)|^p \:dx
&\le &  \int_{\RN}|\gd^l\lp^n u(x) |^p \vpr (x) \: dx \\
&\le & \tilde{C}_0
\end{eqnarray*}
for every $r>0$. By Lebesgue Monotone Convergence Theorem, we conclude that
\begin{equation}
\hskip-.8in
\gd^l\lp^n u \in L^p( \RN) \label {NE-5}
\end{equation}
Now, let us return to the inequality (\ref{NE-2}). Note that
\begin{equation}
\hskip-.8in
\mbox{supp}\{ D^{2m-(2n+l)} \vpr \} \subseteq \{ x\in \RN | R\le |x| \le 2R \}.  \label {supp}
\end{equation}
Then the relations (\ref{NE-2}) and (\ref{supp}) imply
\begin{equation}
\hskip-.8in
\int_{\RN} | \gd^l\lp^n u(x)|^p \vpr(x) \: dx \le \tilde{C}_0^\frac{1}{p'} \left( \int_{R<|x|<2R} |\gd^l\lp^n u(x)|^p \: dx \right)^\frac{1}{p} \label{NE-6}
\end{equation}
where $\tilde{C}_0$ is defined in the manner of (\ref{ctilde}).
Passing to the limit as $R\rightarrow \infty$ in (\ref{NE-6}), by
the absolute convergence of the integral $ \int_{R<|x|<2R}
|\gd^l\lp^n u(x)|^p \:dx$ and the Lebesgue Dominated Convergence
Theorem, we obtain that
\begin{equation*}
\hskip-.8in \int_{\RN}|\gd^l\lp^n u(x)|^p \: dx =0,
\end{equation*}
which concludes our proof.
\end{proof}

\section*{Acknowledgements}
This work was supported in part by the NSF, grants no. DMS-0204794
and DMS-0504619, the BSF grant no. 2004271, the MAOF Fellowship of
the Israeli Council of Higher Education,  the US Civilian Research
and Development Foundation, grant no. RUM1-2654-MO-05, and by the
USA Department of Energy, under contract number W--7405--ENG--36 and
ASCR Program in Applied Mathematical Sciences.

\end{document}